\documentclass{amsart}

\theoremstyle{plain}
\newtheorem{theorem}{Theorem}

\theoremstyle{definition}

\begin{document}

\title[On the cycle structure of hamiltonian regular bipartite graphs]
{On the cycle structure of hamiltonian $k$-regular bipartite graphs of order~$4k$}
\author{Janusz Adamus}
\address{Department of Mathematics, The University of Western Ontario, London, Ontario, N6A 5B7 Canada and
   Institute of Mathematics, Jagiellonian University, Krak{\'o}w, Poland}
\email{jadamus@uwo.ca}

\begin{abstract}
It is shown that a hamiltonian $n/2$-regular bipartite graph $G$ of order $2n>8$ contains a cycle of length $2n-2$. Moreover, if such a cycle can be chosen to omit a pair of adjacent vertices, then $G$ is bipancyclic.
\end{abstract}
\maketitle

In \cite{ES}, Entringer and Schmeichel gave a sufficient condition for a hamiltonian bipartite graph to be bipancyclic.

\begin{theorem}
\label{thm:ES}
A hamiltonian bipartite graph $G$ of order $2n$ and size $\left\|G\right\|>n^2/2$ is bipancyclic (that is, contains cycles of all even lengths up to $2n$). 
\end{theorem}

Interestingly enough, a non-hamiltonian graph with this same bound on the size may contain no long cycles whatsoever. Consider for instance, for $n$ even, a graph obtained from the disjoint union of $H_1=K_{n/2,n/2}$ and $H_2=K_{n/2,n/2}$ by joining a single vertex of $H_1$ with a vertex of $H_2$.\par

In the present note, we are interested in the cycle structure of a hamiltonian bipartite graph of order $2n$, whose every vertex is of degree $n/2$. One immediately verifies that the size of such a graph is precisely $n^2/2$, so the above theorem does not apply. Instead, we prove the following result.

\begin{theorem}
\label{thm:main}
If $G$ is a hamiltonian $n/2$-regular bipartite graph of order $2n>8$, then $G$ contains a cycle $C$ of length $2n-2$. Moreover, if \ $C$ can be chosen to omit a pair of adjacent vertices, then $G$ is bipancyclic.
\end{theorem}

Our motivation for presenting Theorem~\ref{thm:main} is that, although concerning a narrow class of graphs, it plays an important role in the general study of long cycles in balanced bipartite graphs \cite{AA}. We find it also quite amusing that the proof below relies entirely on the combinatorics of the adjacency matrix.

It should be noted that Tian and Zang~\cite{TZ} proved that a hamiltonian bipartite graph of order $2n\geq120$ and minimal degree greater than $\frac{2n}{5}+2$ is necessarily bipancyclic. This result leaves open the case of $|G|=n<60$, in which the above theorem seems to be best to date.

\begin{proof}
Suppose to the contrary that there is a hamiltonian $n/2$-regular bipartite graph $G$ on $2n$ vertices ($n\geq5$), without a cycle of length $2n-2$. Let $X=\{x_1,\dots,x_n\}$ and $Y=\{y_1,\dots,y_n\}$ be the colour classes of $G$, and let $H$ be a Hamilton cycle in $G$; say, $H=x_1y_1x_2y_2\dots x_ny_nx_1$. Let $E=E(G)$ be the edge set of $G$. The requirement that $G$ contain no $C_{2n-2}$ implies that, for every $i=1,\dots,n$,
\[\tag{1}
x_iy_{i-2}\notin E\,,\quad x_iy_{i+1}\notin E\,,\quad\textrm{and}
\]
\[\tag{2}
\textrm{if}\ x_iy_j\in E\ \textrm{for some}\ j\in\{i+2,\dots,n+i-3\}\,,\ \ \textrm{then}\ x_{i+1}y_{j+1}\notin E\,.
\]
(All indices are understood modulo $n$.)\par

Consider the $n\times n$ adjacency matrix $A_G=[a^i_j]_{1\leq i,j\leq n}$, where $a^i_j=1$ if $x_iy_j\in E$, and $a^i_j=-1$ otherwise. Notice that, from adjacency on $H$ and by $(1)$,
\[\tag{3}
a^i_{i-1}=a^i_i=1\quad\textrm{and}\quad a^i_{i-2}=a^i_{i+1}=-1\quad\textrm{for all}\ i,
\]
and by $(2)$,
\[\tag{4}
a^i_j=1\ \Rightarrow\ a^{i+1}_{j+1}=-1\quad\textrm{for}\ i=1,\dots,n,\, j=i+2,\dots,i-3.
\]
As every $x_i$ has precisely $n/2$ neighbours, the entries of each row of $A_G$ sum up to $0$; i.e., $\sum_{j=1}^na^i_j=0$. Therefore, by $(4)$, we also have
\[\tag{5}
a^i_j=-1\ \Rightarrow\ a^{i+1}_{j+1}=1\quad\textrm{for}\ i=1,\dots,n,\, j=i+2,\dots,i-3.
\]
The properties $(3)$, $(4)$ and $(5)$ imply that $A_G$ (and hence $G$ itself) is uniquely determined by the entries $a^1_3,\dots,a^1_{n-2}$, and more importantly, that the sum of entries of the first column of $A_G$ equals
\begin{multline}
\notag
a^1_1+a^1_n+a^1_{n-1}-a^1_{n-2}+a^1_{n-3}-a^1_{n-4}+\dots-a^1_4+a^1_3+a^1_2\\
=a^1_3-a^1_4+\dots+a^1_{n-3}-a^1_{n-2},
\end{multline}
given that $a^1_1+a^1_2+a^1_{n-1}+a^1_n=0$.\par

On the other hand, every column sums up to $0$, as each $y_j$ has precisely $n/2$ neighbours. Hence $\sum_{j=3}^{n-2}a^1_j=0$ and $\sum_{j=3}^{n-2}(-1)^{j+1}a^1_j=0$, and thus $n-4=4l$ for some $l\geq1$, and $\sum_{k=1}^{2l}a^1_{2k+1}=\sum_{k=1}^{2l}a^1_{2k+2}=0$. In general, for any $1\leq i_0\leq n$,
\[\tag{6}
a^{i_0}_{i_0+2}+a^{i_0}_{i_0+4}+\dots+a^{i_0}_{i_0+n-4}=a^{i_0}_{i_0+3}+a^{i_0}_{i_0+5}+\dots+a^{i_0}_{i_0+n-3}=0\,.
\]
Let now $1\leq i_0\leq n$ be such that $a^{i_0}_{i_0+2}=-1$. In fact, we can choose $i_0=1$ or $i_0=2$, for if $a^1_3=1$, then $a^2_4=-1$, by $(4)$. We will show that there exists a $k\in\{3,\dots,n-3\}$ such that
\[
a^{i_0}_{i_0+k}=a^{i_0+k}_{i_0}=1\,.
\]
Suppose otherwise; i.e., suppose that, for all $3\leq k\leq n-3$, $a^{i_0}_{i_0+k}+a^{i_0+k}_{i_0}\in\{0,-2\}$. Notice that, by $(4)$ and $(5)$, $a^{i_0+k}_{i_0}=(-1)^ka^{i_0}_{i_0-k}$ for $k=3,\dots,n-3$. Hence, in particular, $a^{i_0}_{i_0+4}+a^{i_0}_{i_0+n-4}$, $a^{i_0}_{i_0+6}+a^{i_0}_{i_0+n-6},\dots,a^{i_0}_{i_0+2l+2}+a^{i_0}_{i_0+2l+2}$ are all non-positive. In light of $(6)$, this is only possible when $a^{i_0}_{i_0+2}=1$, which contradicts our choice of $i_0$.
\smallskip

To sum up, we have found $i_0$ and $k\in\{3,\dots,n-3\}$ with the property that $a^{i_0-1}_{i_0+1}=a^{i_0}_{i_0+k}=a^{i_0+k}_{i_0}=1$, which is to say that
\[
x_{i_0-1}y_{i_0+1}\in E\,,\ x_{i_0}y_{i_0+k}\in E\,,\ \textrm{and}\ x_{i_0+k}y_{i_0}\in E\,.
\]
Hence a cycle
\[
C\,=\,x_{i_0-1}\,y_{i_0+1}\,x_{i_0+2}\dots y_{i_0+k-1}\,x_{i_0+k}\,y_{i_0}\,x_{i_0}\,y_{i_0+k}\,x_{i_0+k+1}\dots y_{i_0-2}\,x_{i_0-1}
\]
of length $2n-2$ in $G$; a contradiction.
\bigskip

For the proof of the second assertion of the theorem, suppose that $C$ can be chosen so that the omitted vertices $x'$ and $y'$ are adjacent in $G$. Let $G'=G-\{x',y'\}$ be the induced subgraph of $G$ spanned by the vertices of $C$. Then $G'$ is hamiltonian of order $2(n-1)$ and size
\[
\left\|G'\right\|=\left\|G\right\|-(d_G(x')+d_G(y')-1)=n^2/2-n+1\,,
\]
which is greater than $(n-1)^2/2$. Thus $G'$, and hence $G$ itself, is bipancyclic, by Theorem~\ref{thm:ES}.

\end{proof}

%%%%%%%%%%%%%%%%%%%%%%%%%%%%%%%%%%%%%%%%%%%%%%%%%%
\medskip
%%%%%%%%%%%%%%%%%%%%%%%%%%%%%%%%%%%%%%%%%%%%%%%%%%
%References
%%%%%%%%%%%%%%%%%%%%%%%%%%%%%%%%%%%%%%%%%%%%%%%%%%
\bibliographystyle{amsplain}

\end{document}